\documentclass[12pt]{article}
\usepackage{amsfonts}
\usepackage{latexsym}

\def\medskipamount{12pt} \def\smallskipamount{6pt}
\def\arraystretch{1.6}

\newcounter{bitcount}
\newcommand{\bit}[1]{\addtocounter{bitcount}{1}\pagebreak[3]
\subsection{#1}\nopagebreak\setcounter{equation}{0}}

\renewcommand{\theequation}{\thesubsection .\arabic{equation}}
\renewcommand{\thesubsection}{\arabic{bitcount}}

\newcommand{\re}[1]{\mbox{\bf (\ref{#1})}}

\catcode`\@=\active
\catcode`\@=11
\def\@eqnnum{\hbox to .01pt{}\rlap{\bf \hskip -\displaywidth(\theequation)}}
\catcode`\@=12

\begin{document}


\catcode`\@=\active
\catcode`\@=11
\newcommand{\nc}{\newcommand}


\nc{\bs}[1]{ \addvspace{\medskipamount} \pagebreak[3]
\refstepcounter{equation}
\noindent {\bf (\theequation)} \begin{em} \nopagebreak}

\nc{\es}{\end{em} \par \addvspace{\medskipamount} } 

\nc{\ess}{\end{em} \par}

\nc{\br}[1]{ \addvspace{\medskipamount} \pagebreak[3]
\refstepcounter{equation} 
\noindent {\bf (\theequation)} \nopagebreak}

\nc{\brs}[1]{ \pagebreak[3]
\refstepcounter{equation} 
\noindent {\bf (\theequation) #1.} \nopagebreak}

\nc{\er}{\par \addvspace{\medskipamount} }


\nc{\vars}[2]
{{\mathchoice{\mb{#1}}{\mb{#1}}{\mb{#2}}{\mb{#2}}}}
\nc{\C}{\mathbb C}
\nc{\Q}{\mathbb Q}
\nc{\Z}{\mathbb Z}
\renewcommand{\P}{\mathbb P} 


\nc{\oper}[1]{\mathop{\mathchoice{\mbox{\rm #1}}{\mbox{\rm #1}}
{\mbox{\rm \scriptsize #1}}{\mbox{\rm \tiny #1}}}\nolimits}
\nc{\End}{\oper{End}}
\nc{\END}{\oper{\bf End}}
\nc{\Hom}{\oper{Hom}}
\nc{\HOM}{\oper{\bf Hom}}
\nc{\id}{\oper{id}}
\nc{\pdeg}{\oper{pdeg}}
\nc{\Par}{\oper{Par}}
\nc{\Proj}{\oper{Proj}}
\nc{\rk}{\oper{rk}}
\nc{\SPar}{\oper{S\,Par}}
\nc{\Spec}{\oper{Spec}}

\nc{\slantbold}[1]{\mathchoice{{\hspace{-.2ex}\mbox{\boldmath
        $#1$}\hspace{-.2ex}}}{{\hspace{-.2ex}\mbox{\boldmath
        $#1$}\hspace{-.2ex}}}{{\hspace{-.2ex}\mbox{\boldmath
        $\scriptstyle
        #1$}\hspace{-.2ex}}}{{\hspace{-.2ex}\mbox{\boldmath
        $\scriptscriptstyle #1$}\hspace{-.2ex}}}  }
\nc{\E}{\slantbold{E}}
\nc{\F}{\slantbold{F}}
\nc{\R}{\slantbold{R}}
\renewcommand{\H}{\slantbold{H}}
\renewcommand{\d}{\slantbold{d}}
\nc{\h}{\slantbold{h}}

\nc{\D}{\Delta}

\nc{\operlim}[1]{\mathop{\mathchoice{\mbox{\rm #1}}{\mbox{\rm #1}}
{\mbox{\rm \scriptsize #1}}{\mbox{\rm \tiny #1}}}}

\nc{\al}{\alpha}
\nc{\be}{\beta}
\nc{\la}{\lambda}
\nc{\ep}{\varepsilon}
\renewcommand{\i}{\iota}


\nc{\Left}[1]{\hbox{$\left#1\vbox to
    11.5pt{}\right.\nulldelimiterspace=0pt \mathsurround=0pt$}}
\nc{\Right}[1]{\hbox{$\left.\vbox to
    11.5pt{}\right#1\nulldelimiterspace=0pt \mathsurround=0pt$}}


\nc{\beqas}{\begin{eqnarray*}}
\nc{\ci}{{\cal I}}
\nc{\co}{{\cal O}}
\nc{\cx}{{\C^{\times}}}
\nc{\down}{\Big\downarrow}
\nc{\updown}{\Big\updownarrow}
\nc{\downarg}[1]{{\phantom{\scriptstyle #1}\Big\downarrow
    \raisebox{.4ex}{$\scriptstyle #1$}}}
\nc{\eeqas}{\end{eqnarray*}}
\nc{\fp}{\mbox{     $\Box$} \par \addvspace{\smallskipamount}}
\nc{\lrow}{\longrightarrow}
\nc{\pf}{\noindent {\em Proof}}
\nc{\sans}{\, \backslash \,}

\nc{\mod}{\, / \,}

\nc{\M}{{\mathcal M}}
\nc{\N}{{\mathcal N}}
\nc{\A}{{\mathcal S}}

\nc{\beq}{\begin{equation}}
\nc{\eeq}{\end{equation}}

\hyphenation{para-met-riz-ing sub-bundle}

\catcode`\@=12



\noindent
{\LARGE \bf Variation of moduli of parabolic Higgs bundles}
\medskip \\ 
{\bf Michael Thaddeus } \smallskip \\ 
Department of Mathematics, Columbia University \\
2990 Broadway, New York, N.Y. 10027
\renewcommand{\thefootnote}{}
\footnotetext{Partially supported by NSF grant DMS--9808529.} 

\bigskip

A moduli problem in algebraic geometry is the problem of constructing
a space parametrizing all objects of some kind modulo some
equivalence.  If the equivalence is anything but equality, one usually
has to impose some sort of stability condition on the objects
represented.  In many cases, however, this stability condition is not
canonical, but depends on a parameter, which typically varies in a
finite-dimensional rational vector space.  The moduli spaces obtained
for different values of the parameter are birational (at least if
there are any stable points), and for several moduli problems the
birational transformations between the different moduli spaces have
been well characterized.

Without exception, it has been found that the space of parameters
contains a finite number of hyperplanes called {\em walls} on whose
complement the stability condition is locally constant, so that the
moduli space undergoes a birational transformation when a wall is
crossed.  If the moduli spaces are smooth, the birational
transformation typically has the following very special form.  A
subvariety of the moduli space, isomorphic to the total space of a
$\P^m$-bundle, is blown up; the resulting exceptional divisor is a
$\P^m \times \P^n$-bundle over the same base; and it is blown
down along the other ruling to yield the new moduli space, which
therefore contains the total space of a $\P^n$-bundle.

The object of this paper is to extend the narrative above to the
moduli spaces of {\em parabolic Higgs bundles\/} on a curve.  However,
there is a twist in the tale.  The walls, to be sure, still exist and
play their usual role.  But when a wall is crossed, the birational
transformation undergone by the moduli space is {\em not} of the form
described above.  Rather, it is a so-called {\em elementary
  transformation}.  These are birational transformations defined on
varieties admitting a holomorphic symplectic form, in which the
exceptional divisor is a partial flag bundle, indeed a $\P T^*
\P^n$-bundle.  They were discovered by Mukai \cite{m} in the early
1980s, on moduli spaces of sheaves on abelian and K3 surfaces; since
then they have been observed in several settings.  See Huybrechts
\cite{h1,h2} for an informative discussion.  Because of their
symplectic nature, the appearance of these transformations on the
Higgs moduli spaces is quite natural.  It also seems to be related to
the non-triviality of the obstruction space to the moduli problem, a
hint that would be worth pursuing.

The variation of moduli of ordinary parabolic bundles, without a Higgs
field, was studied by Boden and Hu \cite{bh}.  They described the
projective bundles that are the exceptional loci of the birational
transformations.  A paper of the author \cite[\S7]{flip} used
geometric invariant theory to show that the moduli spaces on either
side of the wall become isomorphic after these exceptional loci are
blown up.  The present work, although it describes a similar result
for parabolic Higgs bundles, does not use geometric invariant theory.
Rather, it resembles another work of the author \cite{pair} which
studied similar phenomena in the moduli spaces of so-called Bradlow
pairs.  

More recently, Boden and Yokogawa \cite{by} studied the moduli spaces
of parabolic Higgs bundles which are our present concern, and in
particular computed their Betti numbers.  They found that these are
unchanged by crossing a wall (which also follows directly from our
main result).  This was explained by subsequent work of Nakajima
\cite{n}, who showed that the moduli spaces on either side of a wall
are actually diffeomorphic.  Nakajima's method of proof uses a family
similar to those of Simpson \cite{s1,s2,s3}, containing both of the
moduli spaces in question as fibers.  This is similar to the argument
used by Huybrechts \cite{h1} to prove that {\em compact} holomorphic
symplectic varieties related by an elementary transformation are
diffeomorphic.

Here is an outline of the contents of the paper.  Section 1 reviews a
simple example providing a local model for the elementary
transformations we shall encounter.  Section 2 reviews the definition
of parabolic Higgs bundles, the basic facts about their moduli spaces,
and the chamber structure on the space of parameters, known as {\em
weights}.  Section 3 reviews the deformation theory of parabolic
Higgs bundles.  Section 4 shows how to perform an elementary
modification, in the bundle sense, of a family of parabolic Higgs
bundles along a Cartier divisor, which will be useful in section 6.
Section 5 describes the locus of parabolic Higgs bundles which become
unstable when the weights cross a wall.  This is the locus which must
be removed from the moduli space, and hence the exceptional locus of
the birational map.  Finally, section 6 characterizes the birational
map as the smooth blow-up and blow-down of the moduli space along the
exceptional locus.  The exceptional divisor which dominates this locus
is the partial flag bundle.

A few words on notation.  The boldface letters are reserved for
parabolic Higgs bundles, and for the hyper-objects which arise from
studying them.  Thus $\H$ denotes hypercohomology, $\h$ its dimension,
$\R$ a hyper-direct image, $\HOM$ and $\END$ the two-term complexes
defined in \S3, and so on.  Both parabolic Higgs bundles and two-term
complexes will occasionally be tensored by a line bundle $L$, which means
the obvious thing: that every vector bundle appearing in the
definition gets tensored by $L$.

\bit{An elementary transformation as a pair of geometric quotients}

Let $G = \cx$ act on $X = \C^{m+1} \times \C^{n+1}$ by $\la(u \oplus v) = \la
u \oplus \la^{-1} v$.  If $R = \C[x_0, \dots, x_m, y_0, \dots, y_n]$
is the coordinate ring of $X$, then the action of $G$ induces a
$\Z$-grading $R = \bigoplus_{k \in \Z} R_k$ in which $\deg x_i = 1$
and $\deg y_j = -1$.  The usual affine quotient $X/G$ is then $\Spec
R_0$, which is singular.  Topologically, it is the quotient of $X$ by
the equivalence relation generated by orbit closure.

But there are also two smooth quotients of open subsets of $X$, namely
$X^+ = \Proj \bigoplus_{k \geq 0} R_k$ and $X^- = \Proj \bigoplus_{k
\leq 0} R_k$.  These are both resolutions of $X/G$.  They are
geometric quotients, the former of $(\C^{m+1} \sans 0) \times \C^{n+1}$, and
the latter of $\C^{m+1} \times (\C^{n+1} \sans 0)$, by the $G$-action.
Consequently, $X^+$ is the total space of $\co(-1)^{n+1}$ over $\P^m$,
while $X^-$ is the total space of $\co(-1)^{m+1}$ over $\P^n$.

It is not hard to show that the blow-ups of these bundles along their
zero-sections are isomorphic, with the same exceptional divisor
$\P^m \times \P^n$.  In fact, both are isomorphic to the
fibered product $X^+ \times_{X/G} X^-$: see for example Brion-Procesi
\cite{bp} or the author \cite{flip}.  We shall refer to this whole
construction as the {\em standard example}.

Now suppose that $m=n$.  Inside $X$ there is then a $G$-invariant
hypersurface $Y$ defined by $\sum_{i=0}^n x_i y_i = 0$.  This is the
cone on a smooth quadric, so it is singular only at the origin.
Consequently, its images $Y^+$ and $Y^-$ in the quotients $X^+$ and
$X^-$ are smooth, and the birational map $X^+ \dasharrow X^-$
restricts to one $Y^+ \dasharrow Y^-$.

Since $Y$ contains both axes, $Y^+$ and $Y^-$ contain the
zero-sections which are blown up.  The normal bundle to $Y$,
restricted to $(\C^{n+1} \sans 0) \times 0$, is the trivial bundle given
by evaluating $\sum x_i y_i$.  So the normal bundle to $\P^n$ in
$Y^-$ is the kernel of the tautological map $\co(-1)^{n+1} \to \co$, which
is nothing but $T^*\P^n$.  Likewise the normal bundle to
$\P^n$ in $Y^+$ is also $T^*\P^n$.  The exceptional divisor in
$\tilde{Y}$, the proper transform of $Y^\pm$ in the blow-up of $X$, is
therefore $\P T^* \P^n$, which is the manifold parametrizing
partial flags of type $(1,n)$ in $\C^{n+1}$.  The two natural
projections to $\P^n$ are the restrictions of the blow-downs
$\tilde{Y} \to Y^\pm$ to $\P T^* \P^n$.  

This is the simplest example of an {\em elementary transformation} in
the sense of Mukai \cite{m}.  It differs from the
standard example in that the exceptional divisor is not a product.
Furthermore, the blow-ups are contained in the fibered product $Y^+
\times_{Y/G} Y^-$, but not equal to it.  A dimension count shows that
the product of the exceptional divisors forms another irreducible
component.  

The moduli problems whose variation has been studied in the past
exhibit transformations that locally resemble the standard example
$X^\pm$.  However, as will be seen, the spaces of parabolic Higgs
bundles exhibit transformations that locally resemble the variant
$Y^\pm$.

\bit{Parabolic Higgs bundles}

Now fix, for the remainder of the paper, a smooth complex projective
curve $C$ of genus $g$ with distinguished points $p_1, \dots, p_n$,
where $2g-2+n \geq 0$.
Denote $D$ the effective divisor $p_1 + \cdots + p_n$.

A {\em quasi-parabolic Higgs bundle} $\E$
consists of an algebraic vector bundle $E$ over $C$ equipped with the
following two things.
First, a {\em quasi-parabolic structure} consisting of a full flag
$$0 = E_{i,0} \subset E_{i,1} \subset \cdots \subset E_{i,r} =
E_{p_i}$$ 
in the fiber of $E$ at each $p_i$.  (It is more standard to allow
partial flags, but we do not since our problem is more complicated in
that case.) Second, a {\em Higgs field} $\phi \in H^0 (\End E \otimes
K(D))$ such that $\phi(E_{i,j}) \subset E_{i,j-1} \otimes K(D)$.  A
{\em parabolic Higgs bundle} or PHB is a quasi-parabolic Higgs bundle
equipped further with {\em parabolic weights}, real numbers
\beq \label{a}
1 = \al_{i,0} > \al_{i,1} > \cdots > \al_{i,r} \geq 0
\eeq
associated to each $p_i$.
The {\em parabolic degree} of a PHB is defined to be $\pdeg \E =
\deg E + \sum_{i,j} \al_{i,j}$.  
\er

A subbundle $F$ of a PHB $\E$ can be given a quasi-parabolic structure
simply by intersecting the flags with $F_{p_i}$, and discarding any
subspace $E_{i,j} \cap F_{p_i}$ which coincides with $E_{i,j-1} \cap
F_{p_i}$.  The weights are assigned accordingly.  Likewise, for the
quotient $E/F$, the flags can be projected to $E_{p_i}/F_{p_i}$.  The
weights of $E/F$ are precisely those not assigned to $F$.  A Higgs
field $\phi$ on $E$ also restricts to one on $F$, and projects to one
on $E/F$, provided that $F$ is {\em $\phi$-invariant}, meaning that
$\phi(F) \subset F \otimes K(D)$.  Thus if $F$ is $\phi$-invariant, it
can be regarded as a sub-PHB $\F$.  \er

A PHB $\E$ is said to be {\em semistable} if
for all proper sub-PHBs $\F$, 
\beq \label{b}
\pdeg \F / \rk \F \leq
\pdeg \E / \rk \E,\eeq 
and {\em stable} if the equality is always strict.

For fixed rank $r$, degree $d$, and weights $\al_{i,j}$, 
Yokogawa \cite{y1} 
has constructed a moduli space $\M$ of semistable PHBs
which is a normal quasi-projective variety of dimension
$$ 2r^2(g-1) + 2 + nr(r-1).  $$
As usual in such a moduli problem, it
parametrizes semistable PHBs modulo an equivalence which, on the stable
PHBs, is nothing but isomorphism.  
Standard arguments as in Newstead \cite[\S5.5]{new}
imply the existence of a universal family over the locus of stable PHBs.

For purely numerical reasons, PHBs which are semistable but not stable
can only appear when the weights take special values.  The space of
all possible values for the weights $\al_{i,j}$ can be viewed as a
product $\A = \A_r^n$ of $n$ open simplices of dimension $r$,
determined by \re{a}.  Any PHB which is semistable but not stable with
respect to weights $(\al_{i,j})$ must have an invariant subbundle such
that equality holds in \re{b}.  This means that
$$\rk \F \pdeg \E = \rk \E \pdeg \F$$
and hence 
\beq \label{e} r^+
\Left( d + \sum_{i,j} \al_{i,j} \Right) 
= r \Left( d^+ + \sum_{i,j} n^+_{i,j} \al_{i,j} \Right),
\eeq 
where $r^+ = \rk F$, $d^+ = \deg F$, and the $n^+_{i,j} = \dim
(E_{i,j}\cap F_{p_i}) /(E_{i,j-1} \cap F_{p_i})$, the latter all being
either 0 or 1.  We will refer to $r^+$, $d^+$, and $n^+_{i,j}$ as the
{\em discrete data} $\d^+$ associated to any subbundle of a PHB.  If
the discrete data are fixed, equation \re{e} requires that the point
$(\al_{i,j})$ belongs to the intersection of an affine hyperplane with
$\A$.  Call this intersection a {\em wall}.  Now $r^+$ and the
$n^+_{i,j}$ can only take a finite number of values, since $0 < r^+ <
\rk E$ and $0 \leq n^+_{i,j} \leq 1$, and for each choice of these,
there are only finitely many values of $d^+$ for which the affine
hyperplane touches $\A$.  The walls are therefore finite in number.

\bs{Remark}
\label{f}
The only other discrete data $\d^-$
giving rise to the same wall are $r^- = r - r^+$,
$d^- = d - d^+$, and $n^-_{i,j} = 1 - n^+_{i,j}$.
\es

This is readily verified from \re{e}.  However, it is what fails if
the flags are not full.


On the complement of the walls, the stability condition is evidently
equivalent to semi-stability, and is locally constant.  On each
connected component of the complement, called a {\em chamber}, the
moduli space of stable PHBs is therefore a fixed quasi-projective
variety.  The remainder of the paper is devoted to showing how the
moduli space changes when a wall is crossed.

\bit{Deformation theory of parabolic Higgs bundles}

We will need some basic facts about the deformation theory of PHBs.
This has been worked out very carefully by Yokogawa \cite[2.1, 4.2]{y}
and Markman \cite[6.3]{mark}.  Most of the results we need are special
cases of their work, so we will only sketch the proofs here.

Given bundles $E$, $F$ over $C$ with parabolic structures at $p_1,
\dots, p_n$ and weights $\al_{i,j}$ and $\be_{i,k}$ respectively, a
homomorphism $\rho: E \to F$ is said to be {\em parabolic} if
$\rho(E_{i,j}) \subset F_{i, k-1}$ whenever $\al_{i,j} > \be_{i,k}$,
and {\em strongly parabolic} if $\rho(E_{i,j}) \subset F_{i, k-1}$
whenever $\al_{i,j} \geq \be_{i,k}$.  Let $\Par\Hom(E,F)$ and $\SPar
\Hom(E,F)$ denote the subsheaves of $\Hom(E,F)$ consisting of
parabolic and strongly parabolic homomorphisms, respectively.  Also
let $\Par\End E = \Par \Hom (E,E)$ and $\SPar\End E = \SPar \Hom
(E,E)$.  Note that $\SPar \Hom(E,F(D))$ is naturally dual to
$\Par\Hom(F,E)$ and that any Higgs field $\phi$ belongs
to $H^0(\SPar \End E \otimes K(D))$ by definition.

If $\E$ and $\F$ are PHBs with Higgs fields $\phi$ and $\psi$
respectively, define then a two-term complex $\HOM (\E, \F)$ by
\beq
\label{l}
\Par \Hom (E, F) \lrow \SPar \Hom (E, F) \otimes K(D),
\eeq
with the map given by $ f \mapsto f \phi - \psi f$.  Also let $
\END \E = \HOM (\E,\E)$.  

\bs{Proposition}\label{c} 
The endomorphisms and infinitesimal deformations of a PHB $\E$ are
given by the hypercohomology $\H^0(\END \E)$ and $\H^1(\END\E)$
respectively.  \es

\noindent {\em Sketch of proof}.  Choose trivializations of $E$ on
open sets $V_\al$ covering $C$.   

The local endomorphisms $g_\al \in H^0(V_\al;\Par\End E)$ define a global
endomorphism of $\E$ if and only if first, they agree on the overlaps,
and second, they preserve the Higgs field.  This happens precisely
when the two components of the hypercohomology differential on \v Cech
cochains 
$$C^0(\Par\End E) \lrow C^1(\Par\End E) \oplus C^0(\SPar\End E
\otimes K(D))$$
both vanish on $(g_\al)$.

Let $S = \Spec \C[\ep]/(\ep^2)$ be the spectrum of the ring of dual
numbers.  Then an infinitesimal deformation of $\E$, that is, an
extension of $\E$ to $S \times C$, has transition functions
$1+\ep f_{\al\be}$ and local Higgs fields $\phi + \ep \psi_\al$ for 
$$(f_{\al\be},\psi_\al) \in C^1(\Par\End E) \oplus C^0(\SPar\End E)
\otimes K(D).$$
The compatibility conditions say that this is a \v Cech cocycle, and
equivalent deformations differ by a change of trivialization, which is
a \v Cech coboundary.  \fp

\bs{Proposition}\label{d}
Let $\E$ and $\F$ be stable PHBs such that $\pdeg \E /\rk \E
\geq \pdeg \F / \rk \F$.  Then $\h^0 (\HOM (\E,\F))
= 1$ if $\E$ and $\F$ are isomorphic, and
0 otherwise.
\es

In particular, $\H^0 (\END \E)$ consists only of scalar
multiplications.  \medskip

\noindent {\em Sketch of proof}.  If $\psi \in \H^0(\HOM(\E,\F))$ is
not an isomorphism or zero, then either its kernel or its image
generates a destabilizing subbundle.  If it is an isomorphism
different from a scalar times the identity, then $\psi - \la \id$ is
not an isomorphism or zero for some scalar $\la$.  \fp

\bs{Proposition}
\label{n}
Let $\E$ and $\F$ be families of PHBs over $C$ parametrized by the
same variety $X$, and suppose that $\E_x \cong
\F_x$ for all $x \in X$.  Then there exists a line bundle $L$ over $X$
such that $\F \cong \E \otimes \pi^*L$, where $\pi: C \times X \to X$
is the projection.
\es

\noindent {\em Sketch of proof}.  Let $L$ be the hyper-direct
image $(\R^0 \pi)_* \HOM(\E,\F)$.  This is a line bundle by \re{d},
and it is straightforward to construct the desired isomorphism. \fp

If $\HOM^*$ denotes the complex obtained from $\HOM$ by taking duals and
reversing arrows, then the natural duality mentioned at the beginning
of the section implies
that $\HOM (\E,\F) = \HOM^* (\F,\E) \otimes K$.  By
Serre duality for hypercohomology, it follows that $\H^i (\HOM
(\E,\F)) =  \H^{2-i} (\HOM 
(\F,\E))^*$. Thus if $\E$ is stable, $\h^0(\END \E) = \h^2(
\END \E) = 1$ by \re{d}, so $\h^1(\END \E)$, the dimension of
the deformation space, depends only on the rank and degree of $\E$.
Hence the moduli space is smooth at the stable points, with tangent
space $T_{\E} \M = \H^1 (\END \E)$.

Furthermore, this tangent space is naturally self-dual.  This induces a
nondegenerate holomorphic 2-tensor on the moduli space $\M$.  It turns
out to be alternating and closed, and hence a symplectic form, but we
will not need to know this.

Given bundles $E^+$, $E^-$ over $C$ with parabolic structures at $p_1,
\dots, p_n$, an {\em extension} of $E^-$ by $E^+$ is a short exact
sequence 
$$0 \lrow E^+ \lrow E \lrow E^- \lrow 0,$$
where $E$ has parabolic
structure at $p_1, \dots, p_n$, all morphisms are parabolic, and the
weights of $E$ at $p_i$ are those of $E^+$ together with those of
$E^-$, ordered so that they are decreasing.  (To prevent repeats,
assume that $\al_{i,j} \neq \be_{i,k}$ for all $i,j,k$.) 
A {\em direct sum} of parabolic bundles is, of course, a split
extension.  
An {\em
extension} of PHBs $\E^-$ by $\E^+$ is an extension of the
underlying parabolic bundles as above, together with a Higgs field
$\phi$ on $E$ which restricts to the given Higgs field on $E^+$ (in
particular, $E^+$ is $\phi$-invariant) and projects to the given
Higgs field on $E^-$ (this projection being
well-defined thanks to the $\phi$-invariance of $E^+$).

\bs{Proposition} 
\label{h}
If $\E^+$ and $\E^-$ are PHBs, the extensions of $\E^-$
by $\E^+$, modulo equivalence, are classified by $\H^1 
(\HOM (\E^-,\E^+))$.
\es

\noindent {\em Sketch of proof}.  In terms of local trivializations of
$E^+$ and $E^-$ on an open cover $V_\al$, any extension must have
transition functions and Higgs fields of the forms
\def\arraystretch{1}
$$\left( \begin{array}{cc} 
1 & f_{\al\be} \\  0 & 1
\end{array} \right) \mbox{     and     } 
\left( \begin{array}{cc} 
1 & \psi_\al \\ 0  & 1
\end{array} \right),$$
\def\arraystretch{1.6}
respectively, for a \v Cech cochain
$$(f_{\al\be},\psi_\al) \in 
C^1(\Par\Hom(E^-,E^+)) \oplus C^0(\SPar\Hom(E^-,E^+)\otimes K(D)).$$
The compatibility conditions say that this is a \v Cech cocycle, and
equivalent extensions differ by a \v Cech coboundary.  \fp

For any extension $\E$ of PHBs $\E^-$ by $\E^+$, let 
$\Par\End' E$ and $\SPar\End' E$
be the subsheaves of $\Par\End E$ and $\SPar\End E$ preserving $E^+$,
and let $\END' \E$ be the complex defined as in \re{l}, but using
these subsheaves.  

\bs{Proposition}
\label{m}
The endomorphisms and infinitesimal deformations of the extension
(that is, of $\E$ together with the sub-PHB $\E^+$\hspace{-.4ex}) 
are given by $\H^0(\END' \E)$ and
$\H^1(\END'\E)$ respectively.
\es

\noindent {\em Sketch of proof}.  Similar to that of \re{c}, except that
both the transition functions and the local Higgs fields must be
upper-triangular with respect to parabolic local splittings
$E|_{V_\al} = E^+|_{V_\al} \oplus E^-|_{V_\al}$. \fp

\bit{Elementary modification of families of parabolic Higgs bundles}

Let $\E$ be a family of PHBs over $C$ parametrized by a base scheme
$Y$ of finite type, and let $\i:Z \to Y$ be a Cartier divisor.  Let
$\E^+ \subset \E|_{C \times Z}$ be a family on $Z$ of sub-PHBs, and
let $\E^-$ be the family of quotients.

Let $E'$ be the kernel of the natural surjection $E \to \i_*E^-$ of
coherent sheaves.  As the elementary modification of a locally free
sheaf along a Cartier divisor, this is locally free \cite[2.16]{bbb}.

At every point $x \in Z$, by \re{l} and \re{m} there are deformation maps
$T_xY \to \H^1(\END \E)$ and $T_xZ \to \H^1(\END' \E)$.  Since there
is a short exact sequence of two-term complexes
$$0 \lrow \END'\E \lrow \END \E \lrow \HOM(\E^+,\E^-) \lrow 0,$$
this determines a well-defined map from the normal space to
$\H^1(\HOM(\E^+,\E^-))$.  Since the normal space to a Cartier divisor
is 1-dimensional, this gives a class $\rho_x \in
\H^1(\HOM(\E^+,\E^-))$, well-defined up to a scalar.

\bs{Proposition} 
\label{q}
There is a natural family $\E'$ of PHBs with underlying bundle $E'$,
such that for $x \not\in Z$, $\E'_x = \E^{\phantom{.}}_x$, and for $x \in
Z$, $\E'_x$ is an extension of $\E^+_x$ by $\E^-_x$ with extension
class $\rho_x$.  \es

\pf.  Let $E_{i,0} \subset E_{i,1} \subset \cdots \subset E_{i,r}$ be
the filtration of $E|_{\{ p_i \} \times Y}$ defining the
parabolic structure.  Then the image of each $E_{i,j}$ in
$\i_*(E^-|_{\{ p_i \} \times Z})$ is a subbundle.  Hence $E'_{i,j} =
E_{i,j} \cap E'|_{\{ p_i \} \times Y}$ is an elementary modification
of $E_{i,j}$, and hence a subbundle of $E'|_{\{ p_i \} \times Y}$.
Since these subbundles are nested, this provides a family of parabolic
bundles with underlying bundle $E'$.

Since $E^+$ is a family of sub-PHBs of $E|_{C \times Z}$, the Higgs
field $\phi$ of $\E$ preserves $E'$.  It therefore induces a section
$\phi'$ of $\End E' \otimes K(D)$.  This agrees with $\phi$ away from
$Z$, so there it certainly satisfies the condition $\phi'(E'_{i,j})
\subset E_{i,j-1} \otimes K(D)$.  But this is a closed condition, so
it is satisfied everywhere and $\phi'$ is a Higgs field.  This
provides the desired family $\E'$ of PHBs, which agrees with
$\E$ away from $Z$.

It remains only to prove the last statement.  We claim that the
general case follows from the special case where $Y$ is $S = \Spec
\C[\ep]/(\ep^2)$, the spectrum of the ring of dual numbers, and $Z$ is
the closed point.  For in the general case, given any 
$x \in Z$, choose an embedding $S \to X$ which takes the closed point
to $x$, but whose image is not contained in $Z$.  The
pull-back of the exact sequence 
$$0 \lrow E' \lrow E \lrow \i_*E^- \lrow 0$$
to $C \times S$ is then still exact: since tensoring with the ring of
dual numbers is exact on the right, the only doubtful thing is the
injectivity of the first map, and this follows since $S$ not contained
in $Z$ implies that $E'_S$ and $E^{\phantom{'}}_S$ are isomorphic at
the generic point, and any map of locally free sheaves which is an
isomorphism at the generic point is injective.  And the constructions
of the previous paragraphs for a parabolic structure and Higgs field
on $E'$ clearly commute with this pull-back.

Choose an open cover $V_\al$ of $C$ where $E_x$ splits as a direct sum
of parabolic bundles, 
$E_x|_{V_\al} = E^+_x|_{V_\al} \oplus E^-_x|_{V_\al}$, and extend
this splitting over $S \times V_\al$.  Relative to this splitting,
$E'_x|_{V_\al} = E^+_x|_{V_\al} 
\oplus (E^- \otimes \ci_{(\ep)})_x|_{V_\al}$.  
Since $\E_x$ is an extension, 
the transition functions on $V_\al \cap V_\be$ and the
Higgs field on $V_\al$ respectively have the forms
\def\arraystretch{1}
$$\left(
  \begin{array}{cc} 1 + \ep * & * \\ \ep f_{\al\be} & 1 + \ep *
\end{array} \right) 
\mbox{   and   }
\left( \begin{array}{cc} 
\phi^+ + \ep * & * \\ \ep \psi_\al & \phi^- + \ep *
\end{array} \right),$$ 
where $\phi^\pm$ are the Higgs fields on $E^\pm$, and
$$(f_{\al\be},\psi_\al) \in C^0(\Par\Hom(E_x^+,E_x^-)) \oplus C^1
(\SPar\Hom(E_x^+,E_x^-) \otimes K(D))$$
is a \v Cech cochain representing the class $\rho_x$. 

The fibers of $E'_x$ are spanned by nonvanishing sections of
$E^+$, and by $\ep$ times nonvanishing sections of $E^-$.  Hence to
find the transition functions and local Higgs fields for $\E'_x$, we
must take the coefficients of 
$$\left( \begin{array}{cc} 
\ep^0 & \ep^{-1} \\ \ep^1 & \ep^0
\end{array} \right)$$
in the above matrices.  These are
$$\left( \begin{array}{cc} 
1 & 0 \\ f_{\al\be} & 1
\end{array} \right) \mbox{     and     } 
\left( \begin{array}{cc} 
1 & 0 \\ \psi_\al & 1
\end{array} \right)$$
respectively.  Consequently, $\E'_x$ is precisely the extension of
$\E^+_x$ by $\E^-_x$ with extension class $\rho_x$, as desired.  \fp
\def\arraystretch{1.6}

\bit{What appears and disappears when a wall is crossed}

Choose a point in $\A_r^n$ lying on only one wall $W$.
A small neighborhood of this point touches exactly two chambers, say
$\D^+$ and $\D^-$, and our goal in this section is to see what PHBs
become unstable as we pass from one to the other.  

Say a quasi-PHB is {\em $\D^+$-stable} (resp.\ {\em $\D^-$-stable}) if
it is stable with respect to weights $(\al_{i,j}) \in \D^+$ (resp.\ 
$\D^-$).  Then we seek those $\E$ that are $\D^-$-stable but
$\D^+$-unstable.

To put it another way, let $\M^+$ and $\M^-$ denote the moduli spaces
of $\D^+$-stable and $\D^-$-stable PHBs, respectively.  Since, as we
will see, a generic $\D^-$-stable PHB is also $\D^+$-stable, there is
a birational map $\M^- \dasharrow \M^+$.  We seek to understand the
exceptional locus of this map.

Choosing the wall $W$ is equivalent to choosing
discrete data $\d^+$ so that equality holds in
\re{e} when $(\al_{i,j}) \in W$.  According to $\re{f}$, the
only ambiguity is the possibility of exchanging $\d^+$ for $\d^-$.  By
making this exchange if 
necessary, we may assume without loss of generality that $\geq$ holds
in \re{e} when $(\al_{i,j}) \in \D^+$.  

\bs{Proposition}
\label{k}
If $\E$ is $\D^-$-stable but $\D^+$-unstable, then any
destabilizing bundle has discrete data $\d^+$.  
\es

\pf.  As the weights $(\al_{i,j})$ cross from $\D^+$ to $\D^-$, the
destabilizing subbundle $\E^+$ of $\E$ must cease to destabilize.  Hence
the inequality \re{b} must change from $>$ to $<$.  By \re{f}, this
implies that $\E^+$ has discrete data $\d^+$. \fp

\bs{Proposition}\label{g}
If $\E$ is $\D^-$-stable but $\D^+$-unstable with destabilizing
subbundle $\E^+$, then $\E^+$ and the quotient
$\E^-$ are both $\D^+$-stable and $\D^-$-stable.
\es

\pf.  A $\D^+$-destabilizing subbundle of $\E^+$ would also
destabilize $\E$ and have different discrete data.  Similarly,
the inverse image in $\E$ of a destabilizing subbundle of
$\E^-$ would destabilize $\E$ and have different
discrete data.  Hence $\E^\pm$ are $\D^+$-stable.  If they
were $\D^-$-unstable, the destabilizing subbundle would again have
different discrete data from $\d^+$ or $\d^-$ and so would
violate \re{f}.  \fp

\bs{Proposition}
\label{j}
If $\E$ is $\D^-$-stable but $\D^+$-unstable, then the
$\D^+$-destabilizing sub\-bundle $\E^+$ is unique.  \es

\pf.  Let $\E^-$ be the quotient of $\E$ by $\E^+$.  If $\F$ is
another $\D^+$-destabilizing subbundle, by \re{g} it is
$\D^+$-stable, and by \re{k} it must have the same discrete data, and
in particular, rank.  There is then a nontrivial homomorphism $\F \to
\E^-$ of PHBs, and hence a nontrivial element of $\H^0 (\HOM
(\F,\E^-))$, contradicting \re{d}. \fp

\bs{Proposition}
\label{i}
Let $\E^+$ and $\E^-$ be $\D^+$-stable PHBs with
discrete data $\d^+$ and $\d^-$.  Then any extension of $\E^-$ by $\E^+$
is $\D^+$-unstable, and is $\D^-$-stable if and only if it is not split.
\es

\pf.  The $\D^+$-instability is obvious: $\E^+$ is the destabilizing
sub-PHB.  It is equally clear that if $\E$ splits as $\E^+ \oplus
\E^-$, then $\E^-$ is $\D^-$-destabilizing.  

If the extension $\E$ is $\D^-$-unstable, on the other hand, then
by \re{j} the $\D^-$-destabilizing subbundle $\F$ must not be
$\D^+$-destabilizing, and hence has
discrete data $\d^-$.  The composite map $\F \to \E \to \E^-$ must be a
nontrivial homomorphism of PHBs, since $\F$ and $\E^-$ have the same
incidences with the flags.  Hence there is a nonzero element of
$\H^0 (\HOM (\F, \E^-))$, which by \re{d}
must be an isomorphism, and hence splits $\E$.  \fp

Putting together the last four results, we find:

\bs{Corollary}
\label{o}
If $\E$ is $\D^-$-stable but $\D^+$-unstable, then it can be
expressed uniquely as a nonsplit extension of PHBs
$$0 \lrow \E^+ \lrow \E \lrow \E^- \lrow 0$$
where $\E^\pm$ are stable with discrete data $\d^\pm$.
Conversely, any such extension is $\D^-$-stable but $\D^+$-unstable.
\es

So let $\N^+$ and $\N^-$ be the moduli spaces of semistable (hence
stable) PHBs with discrete data $\d^+$ and $\d^-$, respectively, and
let $\E^\pm$ be universal PHBs over $\N^\pm \times C$.  Define $U^\pm$
to be the hyper-direct images $(\R^1 \pi)_* \HOM (\E^\pm, \E^\mp)$,
where $\pi$ is the projection on $\N^+ \times \N^-$.  By \re{d} and
Serre duality for hypercohomology, $(\R^0 \pi)_*$ and $(\R^2 \pi)_*$
vanish, so $U^\pm$ are vector bundles and are dual to each other.  The
projectivization $\P U^-$ parametrizes all non-split extensions of
PHBs in $\N^-$ by those in $\N^+$.  Applying the Leray spectral
sequence first to the projection $\N^+ \times \N^- \times C \to \N^+
\times \N^-$, then to the projection $\P U^- \to \N^+ \times \N^-$,
gives natural isomorphisms
\beqas 
H^0(\N^+ \times \N^-; \: U^+ \otimes U^-) 
& = & \H^1(\N^+ \times \N^- \times C;
\: U^+ \otimes \HOM (\E^-,\E^+)) \\
& = & \H^1(\P U^- \times C; \: \HOM (\E^-, \E^+) \otimes \co(1)).
\eeqas 
The image of the identity endomorphism has the tautological property
that, for any $x \in \P U^-$, its restriction to $\{ x \} \times C$ is
the class in $\H^1(C;\HOM (\E^-_x, \E^+_x))$ determined by $x$, up to
a scalar.  It
therefore can be used to define a universal extension of $\E^-$ by
$\E^+ \otimes \co(1)$: let
$$(\alpha, \beta) \in C^1(\Par\Hom
(E^-,E^+(1))) \oplus C^0(\SPar \Hom(E^-,E^+(1)) \otimes K(D))$$
be a \v Cech representative; then $\alpha$ defines a family of
extensions of parabolic bundles, and $\beta$ defines a family of Higgs
fields on them, just as in \re{h}.

By \re{i} every PHB in this family is $\d^-$-stable.  So by the
universal property of the moduli space $\M^-$, there exists a morphism
$\P U^- \to \M^-$ whose image is precisely the locus of PHBs which
become unstable when the wall is crossed.  An easy dimension count
shows that $\dim \P U^- < \dim \M^-$, and hence that a generic
$\D^-$-stable PHB is also $\D^+$-stable, as promised.

Let $V^-$ be the cotangent bundle to the fibers of $\P U^-$.  Then
there is certainly a map $\pi^-: \P V^- \to \P U^-$, but the
``Euler sequence'' of the cotangent bundle in this case is
$$ 0 \lrow V^- \stackrel{\pi^+}{\lrow} U^+(-1) \lrow \co \lrow 0,$$
so there is also a map $\pi^+: \P V^- \to \P U^+$.  In fact 
$\P V^-$ is the bundle of
partial flags in $U^-$ of type $(1, \rk U^- - 1)$, and $\pi^\pm$ are
the forgetful morphisms that discard one subspace.

\bs{Proposition}
\label{p}
The morphism $\P U^- \to \M^-$ is an embedding with normal bundle
$V^-$, such that the following diagram commutes:
$$\begin{array}{ccc}
T_x \M^- & \lrow & \H^1 (\END \E_x) \\
\down && \down \\
V^-_x& \stackrel{\pi^+}{\lrow} & \H^1(\HOM(\E_x^+,\E_x^-)).
\end{array}$$

\es

\pf.  It follows from \re{j} that the morphism is injective.  To
show that it is an embedding, we must show that its derivative is also
injective.  

Consider the long exact sequence over $\{ x \} \times C$ of 
$$0 \lrow \HOM(\E^-,\E^+(1)) \lrow \END'\E \lrow \END\E^+ \oplus
\END\E^- \lrow 0.$$
Call these complexes $A^\cdot$, $B^\cdot$, and
$C^\cdot$ respectively.  By Serre duality and \re{d}, $\h^2(C^\cdot) =
2$ and $\h^2(A^\cdot) = 0$.  Therefore $\h^2(B^\cdot)=2$.  On the
other hand, by \re{d} again, $\h^0(B^\cdot)=1$, generated by scalar
multiplications, since $B^\cdot$ is a subcomplex of $\END\E$.  The
connecting homomorphism from $\H^0(C^\cdot)$ to $\H^1(A^\cdot)$
therefore has rank 1.  Its image must be the line spanned by the
extension class $\rho$ of $\E$.  This follows from exactness, since
by \re{m} $\H^1(\END' \E)$ classifies infinitesimal deformations of
extensions, and the deformation of any extension along its extension
class is certainly isomorphic to a trivial one.

On the other hand, $\P U^-$ parametrizes a family of extensions, so
for each $x \in \P U^-$, there is a natural map $T_x \P U^- \to
\H^1(\END'\E_x)$. 
By the previous paragraph this extends to a short
exact sequence of maps:
$$
\begin{array}{ccccccccc}
0 & \lrow & (V^-)^* & \lrow & T\P U^- & \lrow 
& T\N^+ \oplus T \N^- & \lrow & 0 \\
\down && \down && \down && \down && \down \\
0 & \lrow & \H^1(A^\cdot)/\langle \rho \rangle &
\lrow & \H^1(B^\cdot) & \lrow & \H^1(C^\cdot) & \lrow & \phantom{.}0.
\end{array}
$$
The outer maps are clearly isomorphisms, hence so is the middle one by
the 5-lemma.

Now consider the long exact sequence over $\{ x \} \times C$ of 
$$0 \lrow \END'\E \lrow \END\E 
\lrow \HOM(\E^+(1),\E^-) \lrow 0.$$
It shows that $\H^1(\END'\E)$ injects in $\H^1(\END\E)$;
since the former is $T\P U^-$, and the latter is $T\M$, this completes
the proof that the derivative of our morphism is injective.
Furthermore, since $\h^2(\HOM(\E^+(1),\E^-)) = 0$, 
$\h^2(\END\E) = 1$, and $\h^2(\END'\E) = 2$ as seen above, the
connecting homomorphism $\H^1(\HOM(\E^+(1),\E^-)) \to
\H^2(\END'\E)$ has rank 1, and hence the map $\H^1(\END\E)
\to \H^1(\HOM(\E^+(1),\E^-))$, whose image is the normal bundle
we seek, has corank 1.  Now this map is Serre dual to the natural map
$\H^1(\HOM(\E^-,\E^+(1))) \to \H^1(\END\E)$ taking a
deformation of the extension class $\rho$ of $\E$ to a deformation of the
bundle itself.  Since, as observed before, a deformation in the
direction of $\rho$ itself is isomorphic to a trivial deformation, the
kernel of this map contains the line through $\rho$, hence equals it
since it has rank 1.  Therefore the normal bundle is the annihilator
of $\rho$ in $\H^1(\HOM(\E^+(1),\E^-))$, which is exactly $V^-$. \fp

\bit{The elementary transformation of the moduli space}

Let $\tilde{\M}^-$ be the blow-up of $\M^-$ along the image of the
embedding $\P U^- \to \M^-$ of \re{p}.  The exceptional divisor is $\P
V^-$, which is the bundle of partial flags in $U^-$ of type $(1, \rk
U^- - 1)$.  Since $U^+$ is dual to $U^-$, forgetting the 1-dimensional
subspace gives a morphism $\P V^- \to \P U^+$.

\bs{Proposition}
There is a morphism $\tilde{\M}^- \to \M^+$ such that the following
diagram commutes:
$$
\begin{array}{ccccc}
\M^- \sans \P U^- & \lrow & \tilde{\M}^- & \longleftarrow & \P V^- \\
\down && \down && \downarg{\pi^+} \\
\M^+ \sans \P U^+ & \lrow & \M^+ & \longleftarrow & \P U^+.
\end{array}
$$
\es 

\pf.  Let $\E$ be the universal PHB over $\M^-$.  By uniqueness of
families \re{n}, the restriction $\E|_{\P U^- \times C}$ is isomorphic
to the universal extension of $\E^-$ by $\E^+(1)$ defined after
\re{o}, tensored by the pull-back of a line bundle $L$ over $\P U^-$.

So if $\E$ is pulled back to $\tilde{\M}^-$, then its restriction to
the exceptional divisor $\P V^-$ has a family $\E^+ \otimes L(1)$ of
sub-PHBs.

Let $\E'$ be the elementary modification of $\E$ along this family, in
the sense of \re{q}.  Then for $x \not\in \P V^-$, $\E'_x =
\E^{\phantom{.}}_x$, while for $x \in \P V^-$, $\E'_x$ is an extension of
$\E^+_x$ by $\E^-_x$, with extension class given by the image of the
normal space $N_x (\P V^-/\tilde{\M}^-)$ in $\H^1(\HOM(\E^+,\E^-))$ via
the deformation map of $\E$.  This is precisely $\pi^+(x)$, as we see
from the commutative diagram 
$$\begin{array}{ccccc}
T_x \tilde{\M^-} & \lrow & T_{\pi^-(x)} \M^- 
& \lrow & \H^1 (\END \E_x) \\
\down && \down && \down \\
N_x(\P V^-/\tilde{\M^-}) 
& \lrow & V^-_{\pi^-(x)} & \stackrel{\pi^+}{\lrow} 
& \H^1(\HOM(\E_x^+,\E_x^-)).
\end{array}$$

This defines a map $\tilde{\M}^- \to \M^+$ such that the diagram in
the statement commutes; by the universal property of $\M^+$,
it is a morphism. \fp
  
This brings us to our main result.

\bs{Proposition}
There is a natural isomorphism $\tilde{\M}^- \leftrightarrow
\tilde{\M^+}$ such that the following diagram commutes:
$$
\begin{array}{ccccc}
\M^- \sans \P U^- & \lrow & \tilde{\M}^- & \longleftarrow & \P V^- \\
\updown && \updown && \updown \\
\M^+ \sans \P U^+ & \lrow & \tilde{\M}^+ & \longleftarrow & \P V^+.
\end{array}
$$
\es

\pf.  
The roles of the plus and minus signs in the previous proposition
are completely interchangeable, so we have morphisms $\tilde{\M}^- \to
\M^+$ and $\tilde{\M}^+ \to \M^-$.  Combining these with the
blow-downs $\tilde{\M}^\pm \to \M^\pm$ gives injections of both
$\tilde{\M}^+$ and $\tilde{\M}^-$ into $\M^+
\times \M^-$.  Indeed, both injections are embeddings, since as is
easily checked they annihilate no tangent vectors, and both have the
same image, namely the closure of the graph of the isomorphism $\M^-
\sans \P U^- \leftrightarrow \M^+ \sans \P U^+$.  Moreover, the
image of both $\P V^+$ and $\P V^-$ is the incidence correspondence
between $\P U^+$ and $\P U^-$. \fp

The author's previous work \cite{flip} showed that, in the case of
ordinary parabolic bundles, the blow-up of the moduli space $\M^-$ is
actually isomorphic to the fibered product $\M^+ \times_{\M^0} \M^-$.
It is clear from the construction that the blow-up in the case of PHBs
is an irreducible component of the fibered product.  But, as in the
model of \S1, there is another component.  A dimension count shows
that the fibered product of the exceptional divisors has the same
dimension as the blow-up.  

Instead of PHBs with values in the canonical bundle, one could study
several related moduli problems: that of PHBs with values in a bundle
of higher degree, for example, or ``$K(D)$ pairs'' in the sense of
Boden and Yokogawa \cite{by}, or the twistor family of Simpson.  In
those cases the counterparts of $\H^2(\END \E)$ and $\H^2(\END' \E)$
vanish, and hence the normal bundle of the blow-up locus becomes
essentially $\H^1 (\HOM(\E^+(1),\E^-))$.  The story then is more
conventional: the blow-ups and blow-downs locally resemble those of
the standard example, and so on.  The novelty of our situation
therefore seems to be connected to the non-triviality of the
obstruction space of our moduli problem, even though the obstruction
map is of course zero.

A ``master space'' whose quotients, under different linearizations, by
a fixed torus action are the moduli spaces of PHBs with different
weights has not been constructed.  But if there is one, it cannot be
smooth, for then the blow-ups and blow-downs would resemble the
standard example locally, up to a finite cover \cite[1.15]{flip}.
Rather, we might expect it to have ordinary double points, like the
cone on the quadric.


\begin{thebibliography}{99}
{\small

\bibitem{bh}
{\sc H.U. Boden {\rm and} Y. Hu},
Variations of moduli of parabolic bundles, 
{\sl Math.\ Ann.\ }301 (1995) 539--559. 

\bibitem{by} 
{\sc H.U. Boden {\rm and} K. Yokogawa},
Moduli spaces of parabolic Higgs bundles and parabolic $K(D)$ pairs
over smooth curves, I, 
{\sl Internat.\ J. Math.\ }7 (1996) 573--598. 

\bibitem{bp}
{\sc M. Brion {\rm and} C. Procesi}, 
Action d'un tore dans une vari\'et\'e projective, 
{\sl Operator algebras, unitary representations, enveloping algebras, 
and invariant theory}, 
(A. Connes, M. Duflo, A. Joseph, and R. Rentschler, eds.), 
Birkh\"auser (1990)
509--539.

\bibitem{bbb}
{\sc R. Friedman},
{\sl Algebraic surfaces and holomorphic vector bundles,}
Springer, 1998.

\bibitem{gh}
{\sc P. Griffiths {\rm and} J. Harris,} 
{\sl Principles of algebraic geometry,} Wiley, 1978.

\bibitem{h1}
{\sc D. Huybrechts},
Birational symplectic manifolds and their deformations,
{\sl J. Differential Geom.\ }45 (1997) 488--513. 

\bibitem{h2}
{\sc D. Huybrechts},
Compact hyper-K\"ahler manifolds: basic results,
{\sl Invent.\ Math.\ }135 (1999) 63--113.

\bibitem{mark}
{\sc E. Markman}, 
Spectral curves and integrable systems,
{\sl Compositio Math.\ }93 (1994) 255--290.

\bibitem{m}
{\sc S. Mukai},
Symplectic structure of the moduli space of sheaves on an abelian or
$K3$ surface, 
{\sl Invent.\ Math.\ }77 (1984) 101--116. 

\bibitem{n}
{\sc H. Nakajima},
Hyper-K\"ahler structures on moduli spaces of parabolic Higgs bundles
on Riemann surfaces,
{\sl Moduli of vector bundles (Sanda, 1994; Kyoto, 1994)}, 199--208, 
Lecture Notes in Pure and Appl.\ Math.\ 179, Dekker, 1996.

\bibitem{new}
{\sc P.E. Newstead}, 
{\sl Introduction to moduli problems and orbit spaces,} 
Tata Inst., Bombay, 1978. 

\bibitem{s1}
{\sc C.T. Simpson},
Harmonic bundles on noncompact curves, 
{\sl J. Amer.\ Math.\ Soc.\ }3 (1990) 713--770. 

\bibitem{s2}
{\sc C.T. Simpson},
Moduli of representations of the fundamental group of a smooth
projective variety I,
{\sl Inst.\ Hautes \'Etudes Sci.\ Publ.\ Math.\  }79 (1994) 47--129. 

\bibitem{s3}
{\sc C.T. Simpson},
Moduli of representations of the fundamental group of a smooth
projective variety II,
{\sl Inst.\ Hautes \'Etudes Sci.\ Publ.\ Math.\  }80 (1995) 5--79.

\bibitem{pair}
{\sc M. Thaddeus},
Stable pairs, linear systems and the Verlinde formula, 
{\sl Invent.\ Math.\ }117 (1994) 317--353.                  

\bibitem{flip}
{\sc M. Thaddeus},
Geometric invariant theory and flips, 
{\sl J. Amer.\ Math.\ Soc.\ }9 (1996) 691--723. 

\bibitem{y1}
{\sc K. Yokogawa},
Compactification of moduli of parabolic sheaves and moduli of
parabolic Higgs sheaves,  
{\sl J. Math.\ Kyoto Univ.\ }33 (1993) 451--504.

\bibitem{y}
{\sc K. Yokogawa}, 
Infinitesimal deformation of parabolic Higgs sheaves,
{\sl Internat.\ J. Math.\ }6 (1995) 125--148. 

}

\end{thebibliography}
\end{document}